\documentclass[twoside]{amsart}
\usepackage{amsmath}
\usepackage{amssymb}
\usepackage{latexsym}
\usepackage{verbatim}
\usepackage{graphics}

\begin{document}
\newcommand{\R}{{\mathbb R}}
\newcommand{\C}{{\mathbb C}}
\newcommand{\N}{{\mathbb N}}
\newcommand{\T}{{\mathbb T}}
\newcommand{\Si}{{\mathbb S}}
\newcommand{\NN}{{\mathbb N}}
\newcommand{\Z}{\mathbb Z}
\newcommand{\ZZ}{\mathbb Z}
\newcommand{\Q}{{\mathbb Q}}
\newcommand{\A}{{\mathcal A}}
\newcommand{\B}{{\mathcal B}}
\newcommand{\EE}{{\mathcal E}}
\newcommand{\II}{{\mathcal I}}
\newcommand{\QQ}{{\mathcal Q}}
\newcommand{\YY}{{\mathcal Y}}
\newcommand{\XX}{{\mathcal X}}
\newcommand{\OO}{{\mathcal O}}
\newcommand{\HH}{\check{H}}

\numberwithin{equation}{section}

\newtheorem{definition}{Definition}[section]
\newtheorem{theorem}[definition]{Theorem}
\newtheorem{proposition}[definition]{Proposition}
\newtheorem{lemma}[definition]{Lemma}
\newtheorem{corollary}[definition]{Corollary}
\newtheorem{remark}[definition]{Remark}
\newtheorem{conjecture}[definition]{Conjecture}
\newtheorem{example}[definition]{Example}
\newtheorem{problem}[definition]{Open Problem}
\newtheorem{observacion}[definition]{Observation}

\title{Constant-length substitutions and countable scrambled sets}

\author{Fran\c cois Blanchard}

\address{Institut de Math\'ematiques de Luminy (UPR 9016 du CNRS,
FRUMAM)~; case 907, 163 avenue de Luminy, 13288 Marseille Cedex 09,
France}

\email{blanchar@iml.univ-mrs.fr}

\author{Fabien Durand}

\address{ Laboratoire Ami\'enois
de Math\'ematiques Fondamentales  et
Appliqu\'ees, CNRS-UMR 6140, Universit\'{e} de Picardie
Jules Verne, 33 rue Saint Leu, 80000 Amiens, France.}
\email{fdurand@u-picardie.fr}

\author{Alejandro Maass}

\address{ Departamento de Ingenier\'{\i}a
Matem\'atica, Universidad de Chile
and Centro de Modelamiento Ma\-te\-m\'a\-ti\-co,
UMR 2071 UCHILE-CNRS, Casilla 170/3 correo 3,
Santiago, Chile.}
\email{amaass@dim.uchile.cl}

\subjclass{Primary: 58F08; Secondary: 58F03,54H20}

\keywords{substitution, Li--Yorke pairs, topological dynamics}

\begin{abstract}
In this paper we provide examples of topological dynamical systems
having either  finite or countable scrambled sets. In particular we
study conditions for the existence of Li-Yorke, asymptotic and distal
pairs in cons\-tant--length substitution dynamical systems.  Starting
from a circle rotation we also construct a dynamical system having
Li--Yorke pairs, none of which is recurrent.
\end{abstract}

\maketitle

\markboth{ Constant-length substitutions and countable scrambled sets }
{Fran\c cois Blanchard, Fabien Durand, Alejandro Maass}

\section{Introduction}
One definition of topological chaos, based on ideas in  \cite{LY}, emerged
twenty--five years ago. It is not the only definition of chaos by far. It
relies on the existence of uncountable `scrambled' subsets. Here we investigate
what might be called `the edge of Li--Yorke chaos': we give examples of systems
having only finite or countable scrambled subsets. In particular those are
examples of transitive dynamical systems of zero topological entropy 
which are not Li--Yorke chaotic. 

Let $(X,T)$ be a topological dynamical system: $X$ is a compact
metric space  with metric $\varrho$, and $T$ is a surjective
continuous map from  $X$ to itself.  A  pair of points $\{x,y\} \subseteq
X$ is said to be a  {\it Li--Yorke  pair\/} if one has simultaneously
$$
\limsup_{n\to \infty}  \varrho(T^nx,T^ny) >0 \hbox {~~ and~~~ }
\liminf_{n\to \infty}  \varrho(T^nx,T^ny)=0.
$$
A set $S\subseteq X$ is called {\it scrambled\/} if any pair of
distinct  points $\{x,y\}\subseteq S$ is a Li--Yorke pair. Finally, a
system $(X,T)$ is called {\it chaotic in the sense of Li and Yorke\/}
if $X$ contains an uncountable scrambled set. Li--Yorke chaos has 
been recently proved
to result from
various dynamical properties: positive  entropy; 2--scattering;
transitivity together with one periodic orbit  (see \cite{BGKM} and
\cite{HY}). The opposite situation exists too: equicontinuous and
distal systems have no Li--Yorke pairs. The aim of this article is to describe
various systems that are not Li--Yorke chaotic while having Li--Yorke
pairs. Most of them arise from constant-length substitutions, which are the 
most classical topological extensions of odometers.

A pair $\{x,y\} \subseteq X$ is called {\it distal} if  $\liminf_{n\to
\infty} \varrho(T^nx,T^ny)>0$.
If, instead, $\liminf_{n\to
\infty} \varrho(T^nx,T^ny)=0$ then the pair is called {\it proximal}. 
If the limit exists and is equal to zero then the pair is called 
{\it asymptotic}. 
Thus $\{x,y\}$ is a
Li--Yorke pair if and only if it is proximal but not  asymptotic. The
sets of distal, proximal and asymptotic pairs of  $(X,T)$ are denoted
by $\text {\bf D}(X,T)$, $\text {\bf P}(X,T)$ and $\text {\bf
A}(X,T)$ respectively.  Clearly the set
of Li--Yorke pairs is  $\text {\bf LY}(X,T)=\text {\bf P}(X,T)
\setminus \text {\bf  A}(X,T)$.  The sets of distal pairs, Li--Yorke
pairs and asymptotic pairs  partition $X^2$. It is easy to  see that
the image of a proximal (asymptotic) pair under a factor map  is
proximal (asymptotic). Observe that in  all these definitions the
order of a pair is irrelevant: a non--diagonal  pair may be considered
as a subset of cardinality 2. Finally we point out that 
even if definitions of distal, proximal, asymptotic and Li-Yorke pairs work 
for any general dynamical system $(X,T)$, in the examples that will follow we 
essentially look to the case where $T$ is a homeomorphism.

A dynamical system $(X, T)$ is called  {\it minimal} if the unique closed invariant
subsets of $X$ are
$X$ and $\emptyset$.  A {\it distal\/} dynamical system is one in
which every non--diagonal  pair is distal. A dynamical system having
no Li--Yorke pairs is called  {\it almost  distal} \cite{BGKM};
Sturmian systems and the  Morse system are elementary examples. In
\cite{AA} Akin and Auslander introduce {\it
strong} Li--Yorke pairs,  i.e., those Li--Yorke pairs that are  recurrent under
$T\times T$.  They call
{\it semi--distal} a system  without strong Li--Yorke pairs. Distal,
almost distal and semi--distal systems are minimal when  transitive.
A system described some time ago by Floyd (\cite{Au}, p. 26) was
recently remarked to be semi--distal
but not almost distal (\cite{Y}, \cite{AA}); it is an extension
of an adding machine in which fibers are intervals or singletons. Here we
give several examples with the same property, one comes from a substitution of
constant length and the other is a
bounded--to--one extension of an irrational
rotation.

In Section \ref{subst} we study distal, asymptotic and Li--Yorke pairs in
systems generated by constant--length substitutions; this study turns
out to be complete when the alphabet has size 2. This yields examples
of  systems in which scrambled sets have any  given finite
cardinality.  The following section is devoted to the  construction of
an inverse  limit of substitution systems having the  property that all
scrambled sets are at most countable. Finally in Section 5 we construct a
new semi--distal, not almost distal, system.

\section{Preliminaries about substitutions}

\subsection{General facts about substitutions}   A \emph{substitution} is a
map \( \tau :A\rightarrow A^{+} \), where $A$ is a finite set and \(
A^{+} \)  is the set of finite sequences with values in \( A \);
$\tau $ is a  substitution of {\it constant length $p$}, $p \geq 2$,
if $|\tau
(a)| = p$ for  any $a\in A$, where $|\cdot|$ denotes the length of a
word.

In the sequel $\NN$ stands for the set of non-negative integers.
The sets of one--sided and two--sided infinite
sequences are denoted by \(  A^{\N } \) and \( A^{\Z } \) respectively.
If $K=\N$ or $\Z$, points of $A^K$ are denoted by $(x_i)_{i\in K}$ and given
$i,j \in K$, $i \leq j$, $x(i,j)$ denotes the word (or sub--word) $x_i...x_j$
of $x$; by convention
$x(i,i)=x_i$. For  a word $w \in A^+$
if $|w|=n$ we put $w=w_0...w_{n-1}$ and $w(i,j)=w_i...w_j$ for
$0\leq i \leq j < n$. 
The {\it shift map} $T : A^K \to A^K$ is the continuous map defined by $T((x_i)_{i\in K}) = (x_{i+1})_{i\in K}$.
A {\it subshift} $X$ is a closed $T$-invariant subset of
\(  A^{\N } \) or \( A^{\Z } \); the action we consider on $X$
will always be the restriction of the shift map to $X$ and will also be denoted by $T$.
In this paper we will look exclusively at the case where $K=\ZZ$; in this case $T$ is a homeomorphism.
We say that a word appears in a subshift or is
a sub--word of a subshift if it is a sub--word of a point of the
subshift. A subshift is completely determined
by the list of all words that never occur as sub--words of
its points.

The substitution \(  \tau \) can be naturally extended by
concatenation  to \( A^{+} \), \(  A^{\N } \) and \( A^{\Z } \); for
$x=(x_i)_{i\in \Z} \in A^{\Z }$ the extension is given by
$$
\tau(x)=...\tau(x_{-2})\tau(x_{-1}). \tau(x_0)\tau(x_1)...
$$
where the central dot separates negative and non-negative coordinates
of $\tau(x)$.
A further natural  convention is that  the image of the empty word
$\varepsilon$
is $\varepsilon$.  We say that \( \tau \) is \emph{primitive} if there
is $n \in \N$   such that $a$ appears in $\tau^n(b)$ for every $a,b \in
A$.   The substitution $\tau $ generates a subshift  $X_\tau$  of
\( A^{\Z } \) which is the smallest subshift of $A^\Z$ admitting all words
$\{\tau^n(a):  n\in \N, a\in A \}$; when \( \tau \) is primitive  the
subshift it  generates is minimal. For more details and complements about this subsection we refer the reader to \cite{Qu}.
To end this subsection we recall the following result due to B. Moss\'e \cite{Mo} (see also \cite{MS}).

\begin{theorem}
\label{BMosse}
Let $\tau $ be a primitive substitution. Suppose $X_\tau$ is infinite. 
Then, $\tau : X_\tau \to \tau (X_\tau)$ is a one-to-one continuous map. 
If $\tau$ is of 
constant length $p$, then $T^p \circ \tau=\tau \circ T$ and $\tau(X_\tau)$
is a proper $p$-periodic subset of $X_\tau$, that is, 
$\{T^i(\tau(X_\tau)): i \in  \{0,...,p-1\} \}$ is a clopen (closed and open) 
partition of $X_\tau$.
\end{theorem}

In the case where  $\tau$ is of
constant length $p$ the partition 
$\{T^i(\tau(X_\tau)): i \in  \{0,...,p-1\} \}$ will be called
the {\it fundamental partition} of $X_\tau$.

\medskip

Let $\tau$ be a primitive substitution. 
Let us now recall a way to prove that $X_\tau$ is infinite or not.
J.-J. Pansiot \cite{Pa}, and, T. Harju and M. Linna \cite{HL} proved (in our settings) that it is decidable whether $X_\tau$ is infinite or not.
We need some definitions. 
Let $L(\tau )$ be the set of finite words having an occurrence in some $x\in X_\tau$. We say a word $u\in L(X_\tau)$ is {\it biprolongeable} (with respect to $X_\tau$) if there exist two distinct letters $a$ and $b$ such that $ua$ and $ub$ belong to $L(X_\tau)$.
We say $\tau : A\to A^{+}$ is {\it simplifiable} if there exist an alphabet $B$, $|B| < | A|$, and two morphisms $f : A^* \to B^*$, $g : B^* \to A^*$ such that $\tau = g\circ f$.
The substitution $\tau$ is {\it elementary} if it is not simplifiable. 
Note that when $|A| = 2$, $\tau$ is simplifiable if and only if there exist $u\in A^+$, $n,m\in \NN$ such that $\tau (A) = \{ u^n , u^m \}$.

\begin{proposition}
\label{periodic}
\cite{Pa}
Let $\tau : A \to A^+$ be an elementary primitive substitution. Then, $X_\tau$ is infinite if and only if there exists at least one biprolongeable letter $a\in A$. 
\end{proposition}

If $\tau$ is not elementary then there exist an alphabet $B$, $|B| < | A|$, and two morphisms $f : A^* \to B^*$, $g : B^* \to A^*$ such that $\tau = g\circ f$.
Let $\gamma = f\circ g$. We can remark that $X_\gamma$ is infinite if and only if $X_\tau$ is infinite. 

If $\gamma$ is elementary we apply Proposition \ref{periodic} to know if $X_\tau$ is infinite. 
Otherwise we continue the induction. It will stop because the sequence of the cardinality of the alphabets we produce is decreasing and because when $|A|=1$, $X_\tau$ is finite.

\subsection{One--to--one reduction of substitutions}
\label{reduc}

It is often very convenient to assume that a substitution is
one--to--one. Here we show that this can be done without any serious loss of
generality.

   Let $\tau : A \to A^+$ be a substitution, and let $B \subset A$ be such
that for any $a \in A$ there is a unique $b \in B$ for which
$\tau(a)= \tau(b)$. We define the onto map $\phi: A \to B$ by $\phi(a)
= b$  if $\tau (a) = \tau (b)$.  We call $\bar\tau : B \to B^+$ the
unique substitution satisfying $\bar\tau \circ \phi = \phi \circ \tau
$.  If $\tau$ is primitive then $\bar\tau$ is primitive too. We say
$\bar\tau$ is a {\it reduction} of $\tau $.  The map $\phi$ defines a
topological factor from $(X_\tau,T)$  onto $(X_{\bar\tau },T)$,
also called $\phi$. One checks that $\tau ( \phi (a) ) =  \tau (a) $  for
all $a\in A$; thus $\tau(\phi(x))=\tau(x)$ for all $x\in X_\tau$.

   \begin{proposition} \label{reduction} Let $\tau : A \to A^+$ be a
   primitive substitution. The set $X_\tau$ is finite  if and only if
   $X_{\bar\tau}$ is finite. Moreover, if  $X_{\tau}$ is not finite then
   $X_{\tau}$  is topologically  conjugate to  $X_{\bar\tau }$.
\end{proposition}
\begin{proof} If $X_\tau$ is finite $X_{\bar\tau}$ is finite too.
   On the other hand if $X_{\bar\tau}$ is finite all its points
are  periodic for the shift.  Therefore $\phi(x)$ is periodic and
$\tau(x)=\tau(\phi(x))$ too for $x \in X_\tau$. Since
$X_\tau$ is minimal and  $\tau(x)\in X_\tau$ one concludes that $X_\tau$ is
finite.  Suppose $X_\tau$ is not finite. Since $\phi$ is a factor map
we only need  to prove $\phi$ is one--to--one. If $\phi (x) = \phi (y)$,
$x,y\in X_{\tau}$, then  $ \tau (x) = \tau( \phi (x)) = \tau ( \phi
(y)) = \tau (y)$.  But in this case
$\tau : X_\tau \rightarrow \tau (X_\tau)$ is
one--to--one (Theorem \ref{BMosse}) and $x=y$. Therefore $\phi $ is a
topological conjugacy.
\end{proof}

Remark that if $X_\tau$ is finite then $(X_\tau,T)$ is not
   conjugate to $(X_{\bar\tau},T)$. For example take  $\tau:\{0,1\} \to
   \{0,1\}^+$, $\tau(0)=\tau(1)=01$.   A substitution $\tau: A
   \to A^+$ is said to be one--to--one if for all $a,b \in A$,  $a \not = b$,
   $\tau(a) \not = \tau(b)$; $\tau $ is one--to--one if and
   only if $|B| = |A|$. In this case $\tau$ and
   $\bar\tau$ are the same. If $\tau $ is not one--to--one then $|B| <
   |A|$. By repeating the reduction described above finitely many
   times we obtain a one--to--one substitution $\sigma : C\to
   C^+$ and a factor map $\psi : (X_\tau , T) \to (X_\sigma , T) $,
   which is a conjugacy when $X_\tau$ is not finite. The substitution
   $\sigma $ is called the {\it one--to--one reduction of} $\tau$
   (since it is uniquely defined up to the rename of letters in the alphabet).

Note that when $X_\tau$ is finite, the one-to-one reduction can be defined on an alphabet which is not a singleton. It is the case for 
$\tau: \{0,1\} \to
   \{0,1\}^+$ defined by $\tau (0) = 010$ and $\tau (1) = 101$ where $X_\tau$ 
is an orbit of period two.

\subsection{Constant-length substitutions and odometers} 
\label{pi}

Let $p \geq 2$  be an  integer.  The inverse limit $\OO_p$ of the sequence of
groups $(\Z / p^k \Z: k \in \N)$ endowed with the addition of 1 is
called the $p$--odometer; $\OO_P$ is a 
compact topological ring. 
Odometers are always minimal and
uniquely ergodic.  The $p$--odometer will be denoted by
$(\OO_p, T_0)$. With $[p]=\{0,...,p-1\}$ there is a natural homeomorphism 
$n_p: [p]^\N \to \OO_p$ such that $n_p(\delta)=\sum_{i\geq 0} \delta_i p^i$ for 
$\delta=(\delta_i)_{i\in \N} \in [p]^\N$. The sequence $\delta$ is the 
expansion in base $p$ of number $n_p(\delta)$. 
Through this homeomorphism we will identify $[p]^\N$ with $\OO_p$.
The integers can be identified 
with the integer multiples of the identity, and so, we can think in $\Z$
as a dense subring of $\OO_p$. A number in $\OO_p$ is positive 
(respectively negative) integer if and only if its expansion sequence is 
eventually $0$ (resp. eventually $p-1$). For addition in $\OO_p$ you   
add the digits modulo $p$ but with the carry to the right. The projection to 
$\Z_{p^k}$ is given by the list of the first $k$ digits.

Let $\tau:A \to A^+$ be
a primitive substitution  of constant length  $p$ such that
$X_\tau$ is not finite.  The $p$--odometer
is a factor of the subshift $X_\tau$. The proof of this
fact is  based on the following central  result of the theory of
substitutions.

\begin{theorem}[\cite{De}, \cite{Mo}]
\label{Mosse}
Let $\tau$ be a primitive substitution such that $X_\tau$ is not
finite. Then  for any $x \in X_\tau$ there are a unique  sequence
of points  $(x^{(i)})_{i\in \N} \subseteq X_\tau$ and a unique
sequence of positive  integers $(\delta_i)_{i\in \N} \subseteq [p]$
such that $x^{(0)}=x$ and
$\tau(x^{(i)})=T^{- \delta_{i-1}}(x^{(i-1)})$, $i \ge 1$. Moreover, the
map  $\pi:X_\tau  \to [p]^\N$ defined by $\pi(x)=(\delta_i)_{i\in \N}$ is
continuous.
\end{theorem}
\begin{proof}
We sketch part of the proof. 
The fundamental partition is naturally indexed by $[p]$, that is, 
to $j \in [p]$ we associate $T^j(\tau(X_\tau))$. The point 
$x=x^{(0)}$ lies in a unique member of the partition with index $\delta_0$. Pull
back by $T^{\delta_0}$ to get an element of $\tau(X_\tau)$. Then apply 
$\tau^{-1}$ to obtain the element $x^{(1)}$. Continue inductively. 
\end{proof}

For $i \in \N$ and $x \in X_\tau$ put
$p_i(x)=x^{(i)}(-\delta_i,-1)$ and
$s_i(x)=x^{(i)}(1,p-1-\delta_i)$.  If  $\delta_i=p-1$ we put $s_i(x)$
to be the empty word. Analogously if $\delta_i=0$ let
$p_i(x)$ be the empty word. From Theorem \ref{Mosse}, for all
$i\geq 0$  one has
\begin{equation}
\label{eq:centro}
\tau(x_0^{(i+1)})=p_i(x)x_0^{(i)}s_i(x)
\end{equation}
and
\begin{equation}
\label{eq:conmute}
\begin{array}{ll}
|p_i(Tx)|=0 & \hbox{ if } i < i^*, \\
|p_i(Tx)|=|p_i(x)|+1 & \hbox{ if } i=i^*,  \\
|p_i(Tx)|=|p_i(x)| & \hbox{ if } i>i^*, \\
\end{array}
\end{equation}
where $i^*=\inf\{i \in \N: \delta_i\not = p-1\}$.
Set
${\mathcal I}^-(x)=\{ i \in \N: |p_i(x)|\not = 0 \}$ and
${\mathcal I}^+(x)=\{ i \in \N: |s_i(x)|\not = 0 \}$. By definition
if $\pi(x)=\pi(y)$ then ${\mathcal I}^-(x)={\mathcal I}^-(y)$ and
${\mathcal I}^+(x)={\mathcal I}^+(y)$.

We will frequently look at pairs $(x,y)$ such that $\pi (x) = \pi (y)$.
Two such points have a common expansion in 
base $p$ sequence $(\delta_i)_{i\in \NN}$.
In particular it follows that $|p_i (x)| = |p_i (y)|$, $|s_i (x)| = |s_i (y)|$, $i\in \NN$, and ${\mathcal I}^\pm (x) = {\mathcal I}^\pm (y)$. 
In this case we will let 
${\mathcal I} = {\mathcal I}^+ (x) = {\mathcal I}^+ (y)$ 
or ${\mathcal I} = {\mathcal I}^- (x) = {\mathcal I}^- (y)$. 
While the definitions are symmetric, we are interested in long-term positive time behavior and so our attention will focus on ${\mathcal I}^+$ and $s_i$.

\begin{lemma}
\label{vert-hor}
Let $x \in X_\tau$. Then
$$\tau^i(p_i(x))... \tau^2(p_2(x))\tau(p_1(x)) p_0(x) .   x_0 s_0(x)
\tau(s_1(x))\tau^2(s_2(x))... \tau^i(s_i(x))$$
is the sub--word of $x$ centered in $x_0$ and
$T^{-(\delta_0+ \delta_1 p + \cdots + \delta_i p^i)}
x =  \tau^{i+1} (x^{(i+1)})$ for all $i\geq 0$.
\end{lemma}

\begin{proof}
This follows inductively from the definition of $(x^{(i)})_{i\in \N}$ and
$(\delta_i)_{i\in \N}$.
\end{proof}

Define
$k^-(x)=\max
\{|\tau^i(p_i(x))... \tau^2(p_2(x))\tau(p_1(x)) p_0(x)|:
i\in {\mathcal I}^-(x) \}$
and
$k^+(x)=\max\{ |x_0 s_0(x) \tau(s_1(x))\tau^2(s_2(x))... \tau^i(s_i(x))|:
i \in {\mathcal I}^+(x)\}$.
Observe that $k^-(\cdot)$ and
$k^+(\cdot)$ are constants on fibers of $\pi$. 
Also $k^+(x)$ is finite when ${\mathcal I}^+(x)$ is a finite set which is when 
the expansion in base $p$ of $\pi(x)$ is a negative integer in 
$\OO_p$. Furthermore,
if $|{\mathcal I}^+(x)|<\infty$ there is
$j_0 \in \N$ such that $\delta_j=p-1$ for $j \geq j_0$. Therefore
$|s_j(T^{k^+(x)}(x))|=p-1$ for $j \geq 0$
and

\begin{equation}
\label{eq:iplus}
\begin{array}{ll}
x(0,\infty) &=
x(0,k-1) x(k) s_0(T^{k}(x)) \tau(s_{1}(T^{k}(x)))
   \tau^{2}(s_{2}(T^{k}(x)))... \\
&=x(0,k-1)\lim_{m\to \infty} \tau^m(x^{(m)}_1),
\end{array}
\end{equation}

where $k=k^+(x)$.
Analogously if $|{\mathcal I}^-(x)|<\infty$
for all  $j \geq j_0$, $\delta_j=0$ and

\begin{equation}
\label{eq:iminus}
\begin{array}{ll}
x(-\infty,0) &= ... \tau(p_{1}(T^{-k-1}(x) )) p_0(T^{-k-1}(x)) x(-k-1)
x(-k,0) \\
&=
\lim_{m\to \infty}\tau^m(x^{(m)}_{-1}) x(-k,0),
\end{array}
\end{equation}

where $k=k^-(x)$.
Clearly $|{\mathcal I}^+(x)|$ and $|{\mathcal I}^-(x)|$ cannot be
finite at the same time. 

\medskip

Let $u$ and $v$ be two words of $A^+$. Set
$$
[u.v] = \{ (x_n)_{n\in \ZZ} \in X_{\tau} ; x_{-|u|}x_{-|u| +1}\dots
x_{-1} = u, \ x_0 x_1 \dots x_{|v|-1} = v\}.
$$
This is a particular cylinder set where $u$ occurs just
before the zero coordinate and
$v$ begins immediately afterwards.

\begin{lemma}
\label{inter}
Let $(a_nb_nc_n)_{n\in \NN}$ be a sequence on $A^3$ and $(d_n)_{n\in
   \NN}$ be a sequence of integers such that $0\leq d_n \leq p^n -
   1$. Then
$$
\# \bigcap_{n\in \NN} T^{d_n} \tau^n ([a_n.b_nc_n]) \leq 1.
$$
\end{lemma}

The proof is left to the reader.
\smallskip

\begin{lemma}
\label{factorpi}
The map $\pi:(X_\tau,T)\to
(\OO_p,T_0)$ is  a factor map such that  for every $\delta \in \OO_p$,
$|\pi^{-1}(\{\delta\})| < K$, where  $K=|\{$words of length $3$
of $X_\tau \}|$.
\end{lemma}

\begin{proof}
It is clear by (\ref{eq:conmute})
that $T_0\circ \pi=\pi \circ T$, and $\pi$ is
   continuous by Theorem \ref{Mosse}. Therefore, since $(X_\tau,T)$ is minimal,
   $(\pi(X_\tau),T_0)$ is a minimal system and
   $\pi(X_\tau)=\OO_p$. This proves that $\pi$ defines a factor map.
   Let  $\delta = (\delta_i)_{i\in \NN} \in \OO_p$ and suppose
   $|\pi^{-1} (\{ \delta\} )|\geq K+1$.  By Theorem \ref{Mosse} and Lemma
\ref{vert-hor}, for all
   $x\in \pi^{-1} (\{\delta\})$ there  exists a unique sequence of points
   $(x^{(n)})_{n\geq 1}$ such that  $x = T^{(\delta_0 +p\delta_1 +
   \cdots + p^n \delta_n)} \tau^{n+1}(  x^{(n+1)})$ for all $n\in
   \NN$.  Let  $E= \{ x_1,x_2, \dots , x_{K+1}\} \subset \pi^{-1} (
   \{\delta\} )$. Since $|E| > K$,  for all $n\geq 1$ there exists
   $z_1\ne z_2 \in E$  such that $z_1^{(n)}(-1,1) =
   z_2^{(n)}(-1,1)\in A^3$. Then as $E$  is finite
there are $x \ne y\in E$
   and a strictly increasing sequence of integers
   $(n_i)_{i\in \NN}$  such that $x^{(n_i)}(-1,1) = y^{(n_i)}(-1,1)$
   for all $i\in \NN$.  Consequently $x$ and $y$ belong to
   $$T^{(\delta_0 +p\delta_1 + \cdots + p^{n_i-1}  \delta_{n_i-1})}
   \tau^{n_i} [x_{-1}^{(n_i)}. x_{0}^{(n_i)} x_{1}^{(n_i)}].$$
Lemma \ref{inter} implies
   that $x=y$, which contradicts the choice of $x$ and $y$. Thus
   $|\pi^{-1} ( \{\delta\})|\leq K$.
\end{proof}

\section{Li-Yorke pairs of constant--length substitution systems}
\label{subst}
 
Let $\tau:A \to A^+$ be a primitive substitution of  constant length
$p \geq 2$ such that $X_\tau$ is not finite.
The odometer is a distal system (actually an isometry);
whenever $x,y \in X_\tau$ and $\pi(x) \not = \pi(y)$ one has
$(x,y) \in {\bf D}(X_\tau,T)$,
where $\pi$ is the factor map  defined in Theorem \ref{Mosse}.
This implies that
all Li--Yorke pairs of $(X_\tau,T)$, and consequently all scrambled
sets, are contained  in fibres $\pi^{-1}(\{\delta \})$,
$\delta \in \OO_p$.
Then by Lemma \ref{factorpi} all scrambled sets have finite
cardinality.
For the same reason an asymptotic pair is included in one fiber.

A direct consequence of Proposition \ref{reduction} is the following
corollary.

\begin{corollary}
Let $\tau:A \to A^+$ be a primitive substitution of
constant length $p \geq 2$; let $\sigma$
be its one-to-one reduction and $\psi$ be the corresponding map.
If $X_{\tau}$ is not finite
\begin{enumerate}
\item
$\psi\times \psi(  {\bf A}(X_{\tau},T)  ) = {\bf A}(X_{\sigma},T)$,
\item
$\psi\times \psi( {\bf LY}(X_{\tau},T) ) = {\bf LY}(X_{\sigma},T)$,
\item
$\psi\times \psi(  {\bf D}(X_{\tau},T) ) = {\bf D}(X_{\sigma},T)$.
\item
$\psi\times \psi(  {\bf P}(X_{\tau},T) ) = {\bf P}(X_{\sigma},T)$.
\end{enumerate}
Also, $\psi$ maps $\tau$ to $\sigma$. 
\end{corollary}

Therefore to study asymptotic, proximal, distal and Li--Yorke pairs of
$(X_\tau,T)$ it is  enough to consider the one--to--one reduced
system $(X_\sigma,T)$.

\begin{definition}
\begin{enumerate}
\item
Let $u,u' \in A^+$ be two words of the same length. We say they
have a coincidence if $u_i=u'_i$ for some coordinate
$i \in \{0,...,|u|-1\}$. All $i \in \{0,...,|u|-1\}$ such that
$u_i\not = u'_i$ are called non--coincidences.
\item
Let $\tau:A \to A^+$ be a substitution of constant length $p$. We say
that  $\tau$ has a {\rm coincidence} if there are $a,b \in A$,  $a
\not = b$, such that $\tau(a)$ and  $\tau(b)$ have a coincidence.
We say  it
has {\rm overall coincidences} if the last property holds for any $a,b \in A$.
We say it has {\rm  partial coincidence} if there exist
coincidences but  no overall  coincidences.
\end{enumerate}
\end{definition}

When $A=\{0,1\}$ and $\tau$ has coincidences,  then it has
overall coincidences. 

\begin{proposition}
Let $\tau $ be a primitive substitution and assume 
$X_\tau$ is infinite. If $\tau$ has overall coincidences then each 
pair $\{ x,y \} $, $x\not = y$, with $\pi(x)=\pi(y)$  is proximal.
\end{proposition}
\begin{proof}
Let $\{ x,y \} $, $x\not = y$, be a pair such that $\pi(x)=\pi(y)$. 
Shift so as to have ${\mathcal I}$ infinite. Applying Lemma \ref{vert-hor} we see the lengths of the $s_i$ agree. A coincidence occurs within $\tau ( s_i (x))$ and $\tau (s_i (y))$ whenever $\delta_i \not = p-1$. Thus, there is a common 
block of length $p^{i-1}$ between $\tau^i ( s_i (x))$ and $\tau^i (s_i (y))$.
\end{proof}

In the case of the previous proposition, the Li-Yorke pairs are exactly those pairs in a fiber which are not asymptotic.

It is well known that any infinite subshift has non--trivial asymptotic pairs.
Moreover, it has been proved in \cite{Qu} that the number of different
orbits of asymptotic pairs in the Cartesian product is finite (see also
\cite{HZ}).
This property is also true for any infinite subshift with sub--affine
symbolic complexity. In \cite{BDH} the authors give an upper bound for
the number of different orbits of asymptotic pairs for substitution
subshifts.
In the next
propositions we show how asymptotic and Li--Yorke pairs arise
in our substitution systems.

\medskip

\begin{proposition}
\label{super}
Let $\tau: A \to A^+$ be a one--to--one primitive
substitution of constant length
   $p$ such that $X_{\tau}$ is not finite. 
Let $\{ x,y \}$ be a pair such that $\pi (x)=\pi (y)$ and let ${\mathcal I} = {\mathcal I}^+ (x) = {\mathcal I}^+ (y)$ and $k=k^+ (x) = k^+ (y)$.
Then
\begin{enumerate}
\item
$(x,y) \in {\bf A}(X_\tau,T)$ if and only if either 
${\mathcal I}$ is infinite and
$s_i(x)=s_i(y)$ for all sufficiently large $i$, or, 
${\mathcal I}$ is finite and $s_i(T^{k}(x))=s_i(T^{k}(y))$ for 
all sufficiently large $i$.
\item
Assume that $\tau$ has no coincidences. If there exists $i \in \Z$ such that
$x_i \not = y_i$ and $|{\mathcal I}^+(T^i(x))|=\infty$
then it is a  distal pair. Otherwise the pair is asymptotic.
Anyway ${\bf LY}(X_\tau,T) = \emptyset$.
\end{enumerate}
\end{proposition}
\begin{proof}
Let $\pi(x)=\pi(y)=(\delta_i)_{i \in \N }$ and ${\mathcal I}=\{i_0,i_1,...\}$,
where $i_j<i_{j+1}$ for $j \in \N$.
We recall we have
$|s_i(x)|=|s_i(y)|$ for any  $i\in \N$.

\medskip

(1) First assume that $|{\mathcal I}|=\infty$. Then by Lemma
\ref{vert-hor}
$$
x(0,\infty)= x_0
\tau^{i_0}(s_{i_0}(x))\tau^{i_1}(s_{i_1}(x))...;
y(0,\infty)= y_0
\tau^{i_0}(s_{i_0}(y))\tau^{i_1}(s_{i_1}(y))...
$$
Since $\tau $ is one-to-one, $(x,y) \in {\bf A}(X_\tau,T)$
if and only if $s_i(x) = s_i(y)$ for  any $ i \in
{\mathcal I}$ large enough.
Now assume  that $|{\mathcal I}|< \infty$. In this case
$|{\mathcal I}^+(T^k(x))|=\infty$, then by considering
$(T^k(x),T^k(y))$ instead of $(x,y)$ in previous
arguments we conclude (1).

(2) Assume that $\tau$  has no coincidences. If
$x_i \not = y_i$ for some $i \in \Z$ such that
$|{\mathcal I}^+(T^i(x))|=|{\mathcal I}^+(T^i(y))|=\infty$, then
for any $j \in {\mathcal  I}^+(T^i(x))$ all
corresponding symbols of $s_j(T^i(x))$ and $s_j(T^i(y))$, and consequently of
$\tau^j(s_j(T^i(x)))$ and $\tau^j(s_j(T^i(y)))$, are different, which proves
that $(x,y)$ is distal.
Otherwise there is $i \in \Z$ such that $x_i=y_i$ and
$|{\mathcal I}^+(T^i(x))|=|{\mathcal I}^+(T^i(y))|=\infty$.
Then, $s_j(T^i(x))=s_j(T^i(y))$ for any $j \in \N$ which implies
that $(x,y)$ is asymptotic. This completes the proof of (2).
\end{proof}

Let $A$ be an alphabet. In the sequel we will regard $A^2$ as an alphabet of {\it letter-pairs}. 
We identify a pair of words $(u,v)$, with $|u| = |v|$, with a word of $A^2$ in the obvious way.

The next proposition provides a general criterion for existence
of Li--Yorke pairs.

\begin{proposition}
\label{prop:liyorke}
Let $\tau: A \to A^+$ be a one--to--one primitive substitution of
constant length $p$ such that $X_\tau$ is not finite.
Then the system has Li--Yorke pairs if
and only if there exist  $m\in \N$,
$a,b \in A$, $a \not = b$,  such that $\tau^m(a)=uav$, $\tau^m(b)=u'bv'$
with $|u|=|u'|$, $|v|=|v'|$, $v\not = v'$ and $v$ coincides
in at least one coordinate with $v'$.
\end{proposition}
\begin{proof}
First we prove the  condition is necessary.
Let $(x,y) \in {\bf LY}(X_\tau,T)$. Since $\pi(x)=\pi(y)$, one has
$k^+(x)=k^+(y)=k$  and
${\mathcal I}^+(x)={\mathcal I}^+(y)={\mathcal I}$.
If $|{\mathcal I}|<\infty$ instead of
$(x,y)$ consider the pair $(T^k(x),T^k(y))$ which is also a Li--Yorke
pair and for which $|{\mathcal I}(T^k(x))|=\infty$. So we can assume
$|{\mathcal I}|=\infty$; by Lemma \ref{vert-hor}
$x(0,\infty)=x_0\tau^{i_0}(s_{i_0}(x))...
\tau^{i_j}(s_{i_j}(x))...$ and $y(0,\infty)=y_0\tau^{i_0}(s_{i_0}(y))...
\tau^{i_j}(s_{i_j}(y))...$, where $i_0<i_1<...<i_j<...$ are the elements of
${\mathcal I}$.
Since $(x,y)$ is a Li--Yorke pair
the set ${\mathcal J} \subseteq {\mathcal I}$,
${\mathcal J}=\{ i \in {\mathcal I}\ |\ s_i(x)\not = s_i(y)\}$, is infinite.
This implies, by (\ref{eq:centro}), that for any large enough
$i \in \NN$, $x_0^{(i)}\not = y_0^{(i)}$.
There are two cases.

{\bf Case 1: ${\mathcal I}\setminus {\mathcal J}$ is infinite.}

For any $i \in  {\mathcal I}\setminus {\mathcal J}$,
$s_i(x)=s_i(y)$ so $\tau^m(s_i(x))=\tau^m(s_i(y))$ for any $m \in \N$.
Consider large
integers $0<i<j<k$ such that $i,k \in {\mathcal I} \setminus {\mathcal J}$,
$j \in {\mathcal J}$ and $x_0^{(i)} \not = y_0^{(i)}$. Then,
by using (\ref{eq:centro}) several times, we get
$$\tau^{k-i+1}(x^{(k+1)}_0)=u x_0^{(i)} v \tau^{(j-i)}(s_j(x)) w
\tau^{k-i}(s_k(x))$$
and
$$\tau^{k-i+1}(y^{(k+1)}_0)=u' y_0^{(i)} v' \tau^{(j-i)}(s_j(y)) w'
\tau^{k-i}(s_k(y))$$
with $|u|=|u'|$, $|v|=|v'|$ and $|w|=|w'|$.
Then $v \tau^{(j-i)}(s_j(x)) w \tau^{k-i}(s_k(x))$
and $v' \tau^{(j-i)}(s_j(y)) w' \tau^{k-i}(s_k(y))$ are different
(because $j \in {\mathcal J}$) and
have at least one coincidence (because $k \in {\mathcal I}\setminus
{\mathcal J}$). Since
$\mathcal J$ and ${\mathcal I} \setminus {\mathcal J}$
are infinite there exist
triples $i<j<k$ as before such that
$(x^{(k+1)}_0,y^{(k+1)}_0) = (x^{(i)}_0,y^{(i)}_0)$.
Taking $a = x^{(k+1)}_0$,  $b=y^{(k+1)}_0$ and $m=k-i+1$ one
concludes.

{\bf Case 2: ${\mathcal I}\setminus {\mathcal J}$ is finite.}

One  can assume ${\mathcal J}={\mathcal I}$. Define
for any $(a,b) \in A^2$ the set

$$
{\mathcal T}_\tau(a,b)=\{(c,d)\in A^2:
\exists m \leq |A|^2+1, (c,d) \hbox{ occurs in }
(\tau^m(a),\tau^m(b)) \}.
$$

Since $(x,y)$ is a Li--Yorke pair and $s_i(x)\not = s_i(y)$ for any
$i\in {\mathcal I}$, there is an infinite set ${\mathcal R}\subseteq
{\mathcal I}$ such that for any $i \in {\mathcal R}$ there exists
$0\leq j < |s_i(x)|$ with  ${\mathcal T}_\tau((s_i(x))_j,(s_i(y))_j)$
containing a diagonal pair $(c,c)$. That is, for some integer
$0\leq  m \leq |A|^2+1$, $\tau^m((s_i(x))_j)$ and $\tau^m((s_i(y))_j)$
have a coincidence. Since $s_i(x)\not = s_i(y)$ we conclude for large
$i \in {\mathcal R}$ that
$\tau^{m+1}(x^{(i+1)}_0)$ and $\tau^{m+1}(y^{(i+1)}_0)$
have at least one coincidence and two non--coincidences, with
one coincidence and one non--coincidence after a non--coincidence.
An argument similar to the one used to conclude in Case 1
finishes the proof.

\medskip

Let us now prove the condition is sufficient. Assume there exist  $m\in \N$,
$a,b \in A$, $a \not = b$,  such that $\tau^m(a)=uav$, $\tau^m(b)=u'bv'$,
with $|u|=|u'|$, $|v|=|v'|$, $v\not = v'$ and $v$ coincides
in at least one coordinate with $v'$. Without loss of generality suppose
$m=1$. If $|u|>0$ consider the points of $X_\tau$
$$x=... \tau^i(u) ... \tau^2(u) \tau(u) u. a v \tau (v) \tau^2(v) ... \tau^i(v)
   ...$$
and
$$y=.... \tau^i(u') ... \tau^2(u') \tau(u') u'. a v' \tau (v') \tau^2(v') ...
   \tau^i(v') ...$$
where the zero coordinate is just after the central dot.
Clearly $(x,y)$ is a Li--Yorke pair.
If $|u|=0$ taking powers of $\tau$ we can assume there exist
$c,d \in A$ such that $ca$ and $db$ are sub--words of $X_\tau$ and
$\tau(c)=...c$, $\tau(d)=...d$. In this case consider
the points
$$x=\lim_{k \to \infty}\tau^k(c). a v \tau (v) \tau^2(v) ... \tau^i(v) ...;\
y=\lim_{k \to \infty}\tau^k(d). b v' \tau (v') \tau^2(v') ... \tau^i(v') ...$$
which form a Li--Yorke pair.
\end{proof}

The following result can be deduced from
the proof of Proposition \ref{prop:liyorke}.

\begin{corollary}\label{cor:cor1}
Let $\tau: A \to A^+$ be a one--to--one primitive substitution of
constant length $p$ such that $X_\tau$ is not finite.
Assume $\tau$ has overall coincidences. Let $\{ x,y \} \subset X_\tau$ be such that $\pi (x) = \pi (y)$. Then $(x,y) \in {\bf LY}(X_\tau,T)$
   if and only if for infinitely many
$i \in \N$ we have either $s_i(x)\not = s_i(y)$, or,
$s_i(T^{k^+(x)}(x)) \not = s_i(T^{k^+(y)}(y))$.
\end{corollary}

Here are two general sufficient conditions for the existence
of Li--Yorke pairs. The proof is left to the reader.

\begin{corollary}
Let $\tau: A \to A^+$ be a one--to--one primitive substitution of
constant length $p$ such that $X_\tau$ is not finite.
\begin{enumerate}
\item
If there exists $m \in \N$ such that
for any $a,b \in A$, $\tau^m(a)$ and $\tau^m(b)$ have
coincidences and at least two
non--coincidences, then ${\bf LY}(X_\tau,T) \not = \emptyset$.
\item
Assume there exist $m \in \N$ and $a,b \in A$
such that $\tau^m(a)$ and $\tau^m(b)$
have coincidences and at least two
non--coincidences. If for any $(c,c'), (d,d') \in A^2$,
$c \not = c'$, $d\not = d'$,   there is
$i\in \N$ such that $(c,c')$ is a letter-pair of $(\tau^i(d),\tau^i(d'))$
then the system has Li--Yorke pairs.
\end{enumerate}
\end{corollary}

In part (1) of this corollary one cannot replace
`for any $a,b \in A$' by `there are $a,b \in A$'.
Consider the substitution
$\tau: \{0,1,2,3 \} \to \{0,1,2,3 \}^+$ given by
$\tau(0)=0123$, $\tau(1)=1032$, $\tau(2)=1023$ and $\tau(3)=0132$.
For any $m \in \N$
and $a,b \in \{0,1,2,3\}$, either $\tau^m(a)$ and $\tau^m(b)$ have a common
prefix and after this prefix
they have no coincidences, or they have a suffix in common but
no coincidences before this suffix starts.

Then one proves by contradiction that
$\tau$ does not check the condition of Proposition
\ref{prop:liyorke}, and ${\bf LY}(X_\tau,T) = \emptyset$.
On the other hand $\tau^m(1)$ and
$\tau^m(2)$  have coincidences and at least two non--coincidences so
$1$ and $2$ verify the condition in  part (1) of the previous corollary for
any $m \in \N$.

\medskip
Now we give a necessary and sufficient condition to have uncountably many
Li--Yorke pairs.

\begin{proposition}
\label{prop:liyorke2}
Let $\tau: A \to A^+$ be a one--to--one primitive substitution of
constant length $p$ such that $X_\tau$ is not finite.
Then the set of orbits of  Li--Yorke pairs of $X_\tau$
is either finite or uncountable.
Moreover, the set of Li--Yorke pairs is uncountable if
and only if
\begin{equation}
\label{conduncount}
\begin{array}{l}
\exists m \in \N, \ \exists a,b \in A, \ a \not = b,  \tau^m(a)=uavaw, \
\tau^m(b)=u'bv'bw' \\
\hbox{ with  } |u|=|u'|, \ |v|=|v'|, \ |w|=|w'|, \
vaw \hbox{ and } v'bw' \hbox{ have a coincidence.}
\end{array}
\end{equation}
\end{proposition}
\begin{proof}
First we  prove that \eqref{conduncount} is a  necessary and sufficient
condition to have
uncountably many Li--Yorke pairs.

Assume \eqref{conduncount} holds.
Without loss of generality suppose that $m=1$.
We set
$p_0=u$, $p_1=uav$, $s_0=vaw$, $s_1=w$,
$q_0=u'$, $q_1=u'bv'$, $t_0=v'bw'$ and  $t_1=w'$.
Given a sequence ${\underline {\bf n}}=(n_i)_{i\in \N} \in
\{0,1\}^\N$ which contains infinitely many $0$'s and $1$'s
we define points
$$x({\underline {\bf n}})= ... \tau^2(p_{n_2})\tau(p_{n_1})
p_{n_0} .   a s_{n_0} \tau^1(s_{n_1}) \tau^2(s_{n_2}) ... $$
$$y({\underline {\bf n}})= ... \tau^2(q_{n_2})\tau(q_{n_1})
q_{n_0} .   b t_{n_0} \tau^1(t_{n_1}) \tau^2(t_{n_2}) ... $$
where the central dot separates negative and positive coordinates.
One verifies easily
that $x({\underline {\bf n}})$ and $y({\underline {\bf n}})$
belong to $X_\tau$ and that they form a Li--Yorke pair.
Also if  ${\underline {\bf n}},{\underline {\bf n'}}
\in \{0,1\}^\N$, ${\underline {\bf n}} \not = {\underline {\bf n'}}$, then
$x({\underline {\bf n}}) \not = x({\underline {\bf n'}})$ and
$y({\underline {\bf n}}) \not = y({\underline {\bf n'}})$, because
they do not have the same image in the odometer. Therefore,
$\text {\bf LY}(X_\tau,T)$ is uncountable.

\medskip

Now assume that condition \eqref{conduncount}
is not satisfied: then the number of orbits
of Li--Yorke pairs in
$\text {\bf LY}(X_\tau,T)$ is finite.  To prove this claim we consider
a Li--Yorke pair $(x,y) \in {\text {\bf LY}(X_\tau,T)}$ and show that
the sequences
$(x^{(i)}(-\delta_i,-\delta_i+p-1),y^{(i)}(-\delta_i,-\delta_i+p-1))_{i\in \N}$
and $(\delta_i)_{i\in \N}$ are ultimately periodic and periods can be taken
to be smaller than $|A|^2+1$.

As in the proof of Proposition \ref{prop:liyorke} one can assume that
$|{\mathcal I}^+(x)|= |{\mathcal I}^+(y)|= \infty$.
Therefore, for reasons similar to those in Cases 1 and 2
of the proof of Proposition \ref{prop:liyorke}, there exist
$a,b \in A$, $a \not = b$, such that for infinitely many $i\in \N$
there is $m_i\in \N$ for which
\begin{equation}
\label{cond:itera}
\begin{array}{ll}
(i) &   x^{(i)}_0=x^{(i-m_i)}_0=a, \ y^{(i)}_0=y^{(i-m_i)}_0=b, \\
(ii) &   \tau^{m_i}(x^{(i)}_0)=u_i x^{(i-m_i)}_0 v_i, \
\tau^{m_i}(y^{(i)}_0)=u'_i y^{(i-m_i)}_0 v'_i \ ,  \\
  & \hbox{ with } |u_i|=|u'_i|, \ |v_i|=|v'_i|,  \
v_i \hbox{ has a coincidence with } v'_i \ .
\end{array}
\end{equation}

{\it Claim.} For all $r \geq 1$ the word  $(\tau^r(a),\tau^r(b))$
contains $(a,b)$ as a sub-word at most once.
\begin{proof}[Proof of the claim]
Assume this is not true. Then for some
$r \geq 1$, $\tau^r(a)=paqas$, $\tau^r(b)=p'bq'bs'$, with
$|p|=|p'|$ and $|q|=|q'|$. Therefore for any $i \in \N$
such that condition
\eqref{cond:itera} holds
$ \tau^{r+m_i}(a)=uavaw$, $\tau^{r+m_i}(b)=u'bv'bw'$,
  with $|u|=|u'|$,  $|v|=|v'|$,  $|w|=|w'|$  and
$vaw$ has at least one coincidence with $v'bw'$.
This is a contradiction since
condition \eqref{conduncount} does not hold.
\end{proof}

Let $r$ be the smallest positive integer such that
$(\tau^r(a),\tau^r(b))$ has $(a,b)$ as a sub-word.
It results from the claim that $(a,b)$
appears exactly once in $(\tau^{kr}(a),\tau^{kr}(b))$
for  any $k\geq 1$. Moreover, again by the claim, we have that 
$(a,b)$ is not a sub-word of $(\tau^{kr+l}(a),\tau^{kr+l}(b))$ for 
$k\geq 1$ and
$l\in \{1,...,r-1\}$.
Since there are infinitely many $i \in \N$ such that
condition \eqref{cond:itera} holds, this  implies that the sequences
$(x^{(i)}(-\delta_i,-\delta_i+p-1),y^{(i)}(-\delta_i,-\delta_i+p-1))_{i\in \N}$
and $(\delta_i)_{i\in \N}$ are ultimately periodic, and periods can be taken
to be smaller than $|A|^2+1$. This proves that
the set of Li--Yorke pairs is countable.
Moreover, by
\eqref{eq:iplus} and \eqref{eq:iminus},
$(x,y)$ are in the orbit of a Li-- Yorke pair $(x',y')$
such that the sequences
$(x'^{(i)}(-\delta_i,-\delta_i+p-1))_{i\in\N}$,
$(y'^{(i)}(-\delta_i,-\delta_i+p-1))_{i\in \N}$,
are periodic. To conclude we remark, again by
\eqref{eq:iplus} and \eqref{eq:iminus}, that such a sequence determines
a finite number of points; thus there are finitely many
orbits of Li--Yorke pairs in $X_\tau\times X_\tau$.
\end{proof}

\begin{proposition}
\label{recurrent}
Let $\tau: A \to A^+$ be a one--to--one primitive substitution of
constant length $p$ such that  $X_\tau$ is not finite
and condition \eqref{conduncount} holds. Then $(X_\tau,T)$
has strong Li--Yorke pairs.
\end{proposition}
\begin{proof}
By taking a power of the substitution $\tau$
we can assume that
$\exists a,b \in A$,  $a \not = b$,  $\tau(a)=uavaw$,
$\tau(b)=u'bv'bw'$  with   $|u|=|u'|\not = 0$, $|v|=|v'|$,  $|w|=|w'|\not = 0$,
$vaw$  and  $v'bw'$ have a coincidence.
Then the points
$$
x=...\tau^{2m}(u) \tau^{2m-1}(uav) ... \tau(uav)
u . a v a w \tau(w) ... \tau^{2m-1}(w)\tau^{2m}(vaw)...
$$
$$
y=...\tau^{2m}(u') \tau^{2m-1}(u'bv') ... \tau(u'bv')
u' . b v' b w' \tau(w') ... \tau^{2m-1}(w')\tau^{2m}(v'bw')...
$$
form a Li--Yorke pair which is recurrent. Indeed, for every
$n \in \N$ there is $m\in \N$ such that
$(x(-n,n),y(-n,n))$ is a sub-word of $(\tau^{2m}(a),\tau^{2m}(b))$.
\end{proof}

We deduce the following equivalences.

\begin{corollary}
\label{coroequiv}
Let $\tau: A \to A^+$ be a one--to--one primitive substitution of
constant length $p$ such that  $X_\tau$ is not finite. 
Then the following statements are equivalent.
\begin{enumerate}
\item
$(X_\tau , T)$ has uncountably many Li-Yorke pairs;
\item
$(X_\tau , T)$ has infinitely many Li-Yorke orbits;
\item
$(X_\tau , T)$ has at least one strong Li-Yorke pair;
\item
$(X_\tau , T)$ has uncountably many strong Li-Yorke orbits;
\item
Condition \eqref{conduncount} holds for $\tau $. 
\end{enumerate}
\end{corollary}
\begin{proof}
If $(X_\tau , T)$ has finitely many Li-Yorke orbits then $(X_\tau , T)$ has countably many Li-Yorke pairs. Consequently, (1) implies (2).

Proposition \ref{recurrent} gives that (2) implies (3).

We use a remark of Akin to prove that (3) implies (4).
Let $(x,y)$ be a strong Li-Yorke pair of $(X_\tau,T)$. It is
transitive in the closure  of its
orbit, $K=clos(\{(T\times T)^n (x,y): n\in \Z \})$, 
which intersects the
diagonal. Then the
set of transitive points of $(K, T\times T)$ is a dense $G_\delta$
subset of $K$,
therefore uncountable; it is also a set of Li-Yorke pairs of $(X_\tau,T)$.
Hence, $(X_\tau , T)$ has uncountably many strong Li-Yorke pairs and, consequently, uncountably many strong Li-Yorke orbits.

Statement (4) implies (5) by Proposition \ref{prop:liyorke2}. In the same way (5) implies (1).
\end{proof}

{\bf Examples.}

(1) The substitution
$\tau: \{a,b,c  \}\to \{a,b,c   \}^+$, $\tau(a)=aba$,
$\tau(b)=bca$, $\tau(c)=cca$, has countably many
Li--Yorke pairs. Corollary \ref{coroequiv} says that 
$(X_\tau,T)$ has no strong Li--Yorke pairs. 

(2) By Proposition
\ref{prop:liyorke2} the substitution $\tau: \{a,b,c,d\} \to \{a,b,c,d\}^+$
given by
$\tau(a)=baacd$, $\tau(b)=bbbcd$, $\tau(c)=bcaba$ and
$\tau(d)=bdabd$ has uncountably many Li--Yorke pairs;
by Proposition \ref{recurrent}
$(X_\tau,T)$ has strong Li--Yorke pairs.
We show it also has non--recurrent  Li--Yorke pairs.
Consider the points
$$x=\lim_{m \to \infty} \tau^m(b). \tau^m(c); \
y=\lim_{m \to \infty} \tau^m(b). \tau^m(d).
$$
By Corollary \ref{cor:cor1}, $(x,y)$ is a Li-Yorke pair
of $(X_\tau,T)$. But it is not recurrent since
for every positive integers $i,m$, $i<m$,
$(\tau^i(c),\tau^i(d))$ is not a sub-word of
$(\tau^m(c),\tau^m(d))$.

(3)  Let $A=\{0,1\}$ and $\tau: A \to A^+$ be a primitive
substitution with $X_\tau$ not finite. If condition
\eqref{conduncount} holds then
all Li-Yorke pairs of $(X_\tau,T)$ are recurrent.

\bigskip

When $A=\{0,1\}$ we can characterize asymptotic, distal and
Li--Yorke pairs.

\begin{corollary}
\label{coincidence}
Let $\tau: \{0,1\}\to \{0,1\}^+$ be a one--to--one primitive substitution of
constant length $p$ having a coincidence and such that
$X_\tau$ is not finite.
Consider   $\{x,y\}\subseteq
X_\tau$ such that $\pi(x)=\pi(y)$.
\begin{enumerate}
\item
Assume $|{\mathcal I}^+(x)| = \infty$. If for infinitely many  $i \in
{\mathcal I}^+(x)$,  $s_i(x) \not = s_i(y)$, then  $(x,y)$ is a
Li--Yorke pair, otherwise $(x,y)$ is an asymptotic pair.
\item
Assume $|{\mathcal I}^+(x)| <  \infty$. If $\tau(0)$ and $\tau(1)$
begin by the  same letter then $(x,y)$ is an asymptotic pair, otherwise
$(x,y)$ can be an asymptotic or a Li--Yorke pair but not
   a distal pair.
\end{enumerate}
\end{corollary}
\begin{proof}
It is a  consequence of  Proposition \ref{super} (1) and Corollary
\ref{cor:cor1}.
We only give the proof of (2). Assume $|{\mathcal I}^+(x)| <
\infty$ and put $k(x)=k(y)=k$.
If
$s_i(T^k(x))=s_i(T^k(y))$ for any large $i \in \N$,
by Proposition \ref{super} (1), $(x,y)$ is an asymptotic
pair. In particular this condition holds whenever $\tau(0)$ and 
$\tau(1)$ begins by the same letter.
If  $s_i(T^k(x)) \not =s_i(T^k(y))$ for infinitely many $i
\in \N$, since $\tau$ has overall coincidences, then by Corollary 
\ref{cor:cor1}
one concludes that $(x,y)$ is a
Li--Yorke pair.
\end{proof}

\begin{corollary}\label{cor:cor2}
Let $\tau: \{0,1\}\to \{0,1\}^+$ be a one--to--one
primitive substitution of constant length
$p$ such that $X_\tau$ is not finite. 
Consider   $\{x,y\}\subseteq
X_\tau$ such that $\pi(x)=\pi(y)$.
Then

\begin{enumerate}
\item
If there is no coincidence then $\{x,y\}$ 
is distal whenever there is $i \in \N$ such that
$x_i\not = y_i$ and $|{\mathcal
I}^+(T^i(x))|=|{\mathcal I}^+(T^i(y))|= \infty$. In any other  case
the pair is asymptotic.
\item
If there exists a unique non--coincidence then  $\{ x,y \}$ is asymptotic.
\item There
exist Li--Yorke pairs if and only if there exist coincidences  and at
least two non--coincidences.  Moreover, $(x,y)$ is a Li--Yorke pair
if and only if for infinitely many
$i \in \NN$,
$s_i(x) \not = s_i(y)$ if $|{\mathcal I}^+(x)|=\infty$, or,
$s_i(T^{k^+(x)}(x)) \not = s_i(T^{k^+(y)}(y))$
if $|{\mathcal I}^+(x)|<\infty$.
\end{enumerate}
\end{corollary}
\begin{proof}
\begin{enumerate}
\item
Immediate by Proposition \ref{super} (2) since
the alphabet has size two.

\item 
First
assume that $|{\mathcal I}^+(x)|= |{\mathcal I}^+(y)|= \infty$. If
$x^{(i)}(0)=y^{(i)}(0)$  for any $i \in \N$ then $(x,y)$ is an asymptotic
pair. Otherwise for  any large enough $i \in {\mathcal I}^+(x)$ one has
$x^{(i)}(0)\not =  y^{(i)}(0)$; thus, since $\tau$ has a unique
non--coincidence,  $s_i(x)=s_i(y)$. Then by Proposition \ref{super} (1) $x$
and $y$ are asymptotic. If $|{\mathcal I}^+(x)|= |{\mathcal I}^+(y)|<
\infty$ then $|{\mathcal  I}^+(T^k(x))|= |{\mathcal I}^+(T^k(y))|=
\infty$ where $k^+(x)=k^+(y)=k$.
By the argument above  $T^k(x)$ and $T^k(y)$ are asymptotic, and
consequently $x$ and $y$ too.

\item This is a consequence of Corollary \ref{cor:cor1} and Proposition
\ref{prop:liyorke}.
\end{enumerate}
\end{proof}

{\bf Examples.} The system defined by the substitution $\tau:
\{0,1\}\to \{0,1\}^+$, $\tau(0)=010,\ \tau(1)=100$ has  Li--Yorke pairs
and all scrambled sets are finite; the system defined  by the Morse
substitution $\tau: \{0,1\}\to \{0,1\}^+$, $\tau(0)=01,  \tau(1)=10$
has only distal and asymptotic pairs in the fibres over the odometer;
and the system defined by $\tau:  \{0,1\}\to
\{0,1\}^+$, $\tau(0)=01, \tau(1)=00$ has only asymptotic  pairs in fibres
over the odometer. The
last substitution defines a Toeplitz subshift.

\medskip

By increasing the alphabet it is easy to obtain
substitution systems  with scrambled sets having exactly the
cardinality of the alphabet. For $A=\{0,...,n\}$, $n >0$, consider
the primitive one--to--one substitution $\tau: A \to A^+$
defined for $a\in A$ by
$\tau(a)= 0 a a (a+1) 0$ where $n+1=0$. Clearly $X_\tau$ is not finite.
Define for $a \in A$
$$x(a)= ... \tau^i(0) ... \tau(0) 0. a a (a+1) 0 \tau(a (a+1) 0) ...
\tau^i(a (a+1) 0)... \in X_\tau$$
where the zero coordinate is the one just after the central dot.
Then
$S=\{x(a): a\in A \}$ is a scrambled set of size $|A|$.

\section{A dynamical system with countable scrambled sets}

Let $ n \geq 1 $. For  $ A_n = \{ 0,...,n \}$ define the
primitive substitution
$\tau_n: A_n \to A_n^3$ by $\tau_n(a)=a 0 (a+1)$ if
$a \not = n$ and $\tau_n(n)=n 0 n$.  Set
$X_{\tau_n}=X_n$; the shift map on $X_n$  is denoted by $T_n$.
We verify that $X_n$ is not finite. Let
$\rho_n:A_{n+1} \to A_n$ be defined by $\rho_n(a)=a$ if $a \not
= n+1$ and $\rho_n(n+1)=n$; obviously
$\rho_n  \circ \tau_{n+1}=\tau_{n} \circ \rho_n$.
The map $\rho_n$ extends by concatenation to any word
of $A_{n+1}^+$;
it induces a factor map that we
also call  $\rho_n:X_{n+1} \to X_n$. Remark that each $\tau_n$ satisfies the
conditions of Proposition \ref{prop:liyorke}.

\medskip
Let $(X,T)$ be the minimal dynamical system defined by
$$X=\{ (x_n)_{n\geq 1} \in \prod_{n \geq 1}X_n : \forall n \geq 1, \
\rho_n(x_{n+1})=x_n  \}\ ;$$
put
$T((x_n)_{n\geq 1})=(T_n(x_n))_{n\geq 1}$.
Each system  $(X_n,T_n)$ is a factor of $(X,T)$ and the maximal
equicontinuous factor of $(X,T)$ and $(X_n,T_n)$
is the $3$--odometer (see \cite{Qu}).
Denote by $\pi_n:X_n \to \OO_3$ and
$\pi:X\to \OO_3$ the corresponding factor maps.
If $x=(x_n)_{n\geq 1} \in X$ then
$\pi_n(x_n)= \pi(x)$ for all $n$. 
Each $\tau_n$ has overall coincidences, then $\pi_n$ is a proximal map. The 
inverse limit of proximal maps is a proximal map, so, once again, the 
Li--Yorke pairs are exactly the non-asymptotic pairs in the fibers of $\pi$.
We
show that fibres of $\pi$ are at most countable and construct
a countable scrambled set. Since any scrambled set is included in one fibre,
all scrambled sets of $(X,T)$ are finite or countable.

For any  $z \in X_m$, $m \geq 1$, let
$(z^{(i)})_{i\in \N}$ be the sequence in $X_m$
given by Theorem \ref{Mosse}. Recall that if
$\pi_m(z)=(\delta_i)_{i\in \N}$
one has
\begin{equation}
\label{eq:centro4}
\tau_m(z^{(i+1)}_0)=z^{(i)}(-\delta_i,-\delta_i+2)
\end{equation}
and thus
\begin{equation}
\label{eq:down}
z^{(j)}(-\delta_j,-\delta_j+2), j<i,
\hbox{ are determined once we know } z^{(i)}(-\delta_i,-\delta_i+2).
\end{equation}
Both properties are extensively used in the proofs of this
section.

\begin{lemma}
\label{le:prepreimagen}
Let  $x \in X_n$ for some given $n \geq 1$.
\begin{enumerate}
\item If for infinitely many $i \in \NN$,
$x^{(i)}(-\delta_i,-\delta_i+2)=a_i0(a_i+1)$ for some
$a_i \in A_{n-1}$, then there is a sequence
$v \in (A_{n+1}^3)^\NN$ such that
$(y^{(i)}(-\delta_i,-\delta_i+2))_{i \in \N}=v$ for any
$y \in \rho_n^{-1}(\{ x \})$. Moreover, for infinitely many $i \in \N$,
$v_i=a_i0(a_i+1)$ for some $a_i \in A_{n}$.
\item If for any big enough $i \in \N$,
$x^{(i)}(-\delta_i,-\delta_i+2)=n0n$,
then there are two sequences
$v,v' \in (A_{n+1}^3)^\N$ with
$v_i=(n 0 (n+1))$ and $v'_i=((n+1)0(n+1))$ for any enough large
$i \in \N$,
such that
$(y^{(i)}(-\delta_i,-\delta_i+2))_{i \in \N}=v$ or $v'$
for any $y \in \rho_n^{-1}(\{ x \})$.
Moreover, if for infinitely many $i \in \N$, $\delta_i=2$, then
$(y^{(i)}(-\delta_i,-\delta_i+2))_{i \in \N}=v'$
for any $y \in \rho_n^{-1}(\{ x \})$.
\end{enumerate}
\end{lemma}
\begin{proof}
Let $y \in X_{n+1}$ and $x=\rho_{n}(y)\in X_n$.
Put $\pi_n(x)=\pi_{n+1}(y)=(\delta_i)_{i\in \N}$,
then
$\rho_n(y^{(i)}(-\delta_i,-\delta_i+2))=x^{(i)}(-\delta_i,-\delta_i+2)$
for all $i\in \N$.

(1) If $x^{(i)}(-\delta_i,-\delta_i+2)=a0(a+1)$ with
$a \in A_{n-1}$ then $a0(a+1)$ is the unique sub-word
of $X_{n+1}$ such that $\rho_n(a0(a+1))=a0(a+1)$, then
$y^{(i)}(-\delta_i,-\delta_i+2)=a0(a+1)$.
Therefore, if for infinitely many $i\in \N$,
$x^{(i)}(-\delta_i,-\delta_i+2) = a0(a+1)$ for some
$a \in A_{n-1}$, by (\ref{eq:down}),
the sequence $(y^{(i)}(-\delta_i,-\delta_i+2))_{i\in \N}$ is completely
determined given $x$, it is unique and for infinitely many $i\in \N$,
$y^{(i)}(-\delta_i,-\delta_i+2)= a0(a+1)$ for some
$a \in A_{n-1} \subseteq A_n$. This completes the proof of (1).

(2) Assume that $x^{(i)}(-\delta_i,-\delta_i+2)=n0n$
for any $i \geq i_0$.  As $n$ occurs
only in the first and last positions of $\tau_n(n)$,
so $\delta_i$ can only be $0$ or $2$ for $i>i_0$.
Also, if $x^{(i)}(-\delta_i,-\delta_i+2)=n0n$ then
$n0(n+1)$ and $(n+1)0(n+1)$ are the only sub-words
of $X_{n+1}$ such that $\rho_n(n0(n+1))=\rho_n((n+1)0(n+1))=n0n$, then
for $i\geq i_0$,
$y^{(i)}(-\delta_i,-\delta_i+2)=n0(n+1)$ or $(n+1)0(n+1)$.

Assume that for some
$i \geq i_0$, $y^{(i)}(-\delta_i,-\delta_i+2)=n0(n+1)$. Then,
by definition of $\tau_{n+1}$,
$y^{(i+1)}(-\delta_{i+1},-\delta_{i+1}+2)=n0(n+1)$ and
$\delta_{i+1}=0$. Inductively we deduce that
$y^{(j)}(-\delta_j,-\delta_j+2)=n0(n+1)$ and $\delta_{j}=0$ for
any $j > i$. By using (\ref{eq:down}), we conclude
$(y^{(i)}(-\delta_i,-\delta_i+2))_{i\in \N}$ can be equal to two
possibles sequences. One is such that
$(y^{(i)}(-\delta_i,-\delta_i+2))_{i\geq j_0}=(n0(n+1))_{i\geq j_0}$
and the other verifies
$(y^{(i)}(-\delta_i,-\delta_i+2))_{i\geq j_0}=((n+1)0(n+1))_{i\geq j_0}$,
where $j_0 \in \N$. But, if for infinitely many $i\in \N$,
$\delta_i=2$, the sequence $(y^{(i)}(-\delta_i,-\delta_i+2))_{i\in \N}$
can take a unique value such that
$(y^{(i)}(-\delta_i,-\delta_i+2))_{i\geq j_0}=((n+1)0(n+1))_{i\geq j_0}$
for some $j_0 \in \N$.
\end{proof}

\begin{lemma}
\label{le:preimagen}
Let  $x \in X_n$ for some given $n \geq 1$.
\begin{enumerate}
\item   $\rho_n^{-1}(\{ x \})$ contains at most two
elements; when it contains two, one of them has a unique preimage under
$\rho_{n+1}$.
\item  If $x$ has a
unique preimage under $\rho_n$, there exists a unique sequence
\hfill\break
$(x_m)_{m\geq n}$ such that $x=x_n$, $x_m \in X_m$
and $\rho_m(x_{m+1})=x_m$ for
$m \geq n$.
\end{enumerate}
\end{lemma}
\begin{proof}
Let $y \in X_{n+1}$ and $x=\rho_{n}(y)\in X_n$.
Put $\pi_n(x)=\pi_{n+1}(y)=(\delta_i)_{i\in \N}$.
\medskip

(1) We consider three cases.

(a) First suppose that $\delta_i=0$ for $i\geq i_0$.
By definition of $\tau_n$ and $\tau_{n+1}$,
$x^{(i)}_{-1}=n$ and $y^{(i)}_{-1}=n+1$ for $i\geq i_0$.
Then, we deduce from (\ref{eq:iminus}) and Lemma
\ref{vert-hor} that
$x=\lim_{m\to \infty} \tau_n^m(n)x(M,\infty)$ and
$y=\lim_{m\to \infty} \tau_{n+1}^m(n)y(M,\infty)$
for some $M \in \Z$,
where $x(M,\infty)$, $y(M,\infty)$, are uniquely determined by
$(x^{(i)}(-\delta_i,-\delta_i+2))_{i\in \N}$ and
$(y^{(i)}(-\delta_i,-\delta_i+2))_{i\in \N}$ respectively.
Therefore, by Lemma \ref{le:prepreimagen},
$x$ has at most two preimages under
$\rho_n$, one of which has a unique preimage under $\rho_{n+1}$
since it verifies condition (1) in Lemma \ref{le:prepreimagen}.

(b) When $\delta_i=2$ for $i\geq i_0$, by definition
of $\tau_n$, one has $x^{(i)}(-\delta_i,-\delta_i+2)=n0n$
for $i\geq i_0$. We also deduce from (\ref{eq:iplus}) and Lemma
\ref{vert-hor} that for some $M\in \N$ and $a \in A_n$,
$x=x(-\infty,M)\lim_{m\to \infty} \tau_n^m(a)$. Therefore
$y=y(-\infty,M)\lim_{m\to \infty} \tau_{n+1}^m(a)$ if
$a \not = n$, and
$y=y(-\infty,M)\lim_{m\to \infty} \tau_{n+1}^m(n)$
or
$y=y(-\infty,M)\lim_{m\to \infty} \tau_{n+1}^m(n+1)$ if
$a=n$.
In both cases $y(-\infty,M)$ is uniquely determined from
the sequence $(y^{(i)}(-\delta_i,-\delta_i+2))_{i\in \N}$.
But in the case we are studying, by Lemma \ref{le:prepreimagen} (2),
this sequence is uniquely determined from $x$ and verifies
$(y^{(i)}(-\delta_i,-\delta_i+2))_{i\geq j_0}=((n+1)0(n+1))_{i\geq j_0}$
for some $j_0 \in \N$. We conclude $y$ can take at most
two possible values, which proves the first statement of (1).
To finish the proof of (1) in this case observe that
if
$y=y(-\infty,M)\lim_{m\to \infty} \tau_{n+1}^m(n)$ then
$z=z(-\infty,M)\lim_{m\to\infty}\tau^m_{n+2}(n)$
is the unique preimage of $y$ under $\rho_{n+1}$,
where $z(-\infty,M)$ is uniquely determined given $y$ by
Lemma \ref{le:prepreimagen} (2).

(c) Finally if none of the two
conditions above is checked, $x$ and
$y$ are uniquely determined
by the sequences
$(x^{(i)}(-\delta_i,-\delta_i+2))_{i\in \N}$ and
$(y^{(i)}(-\delta_i,-\delta_i+2))_{i\in \N}$. By
Lemma \ref{le:prepreimagen} (2),
$y$ can take at most two values,
one of which must have a unique preimage
under $\rho_{n+1}$ by the same lemma.
This completes the proof of (1).

\medskip

(2) The proof above shows that $y$ is the unique preimage of $x$
under $\rho_n$ if and only if one of the following three cases occur:

(i) for infinitely many $i\in \N$,
$x^{(i)}(-\delta_i,-\delta_i+2)=y^{(i)}(-\delta_i,-\delta_i+2)=
a0(a+1), a \in A_{n-1}$ and $\delta_i\not = 2$;

(ii) for infinitely many $i\in \N$,
$\delta_i=2$, for infinitely many $i\in \N$,
$\delta_i=0$ and there is $i_0 \in \N$ such that
$x^{(i)}(-\delta_i,-\delta_i+2)=n0n$,
$y^{(i)}(-\delta_i,-\delta_i+2)=(n+1)0(n+1)$ for
$i \geq i_0$;

(iii) there is $i_0 \in \N$ such that for $i \geq i_0$,
$\delta_i=2,$
$x^{(i)}(-\delta_i,-\delta_i+2)=n0n$,
$y^{(i)}(-\delta_i,-\delta_i+2)=(n+1)0(n+1)$
and for some $M\in \N$ and $a \in A_{n-1}$,
$y=y(-\infty,M)\lim_{m\to\infty}\tau^m_{n+1}(a)$.

In case (i), by Lemma \ref{le:prepreimagen} and
the proof of (a) and (c) above,
there is a unique $z \in X_{n+2}$ such that
$\rho_{n+1}(z)=y$; it is uniquely determined by
$(y^{(i)}(-\delta_i,-\delta_i+2))_{i\in \N}$. Moreover,
$(z^{(i)}(-\delta_i,-\delta_i+2))_{i\in \N}$ also verifies
condition (i). Then we conclude by induction.

In cases (ii) and (iii) the proof is analogous.
\end{proof}

\begin{proposition}
The map $\pi:X \to \OO_3$ is at most countable--to--one.
\end{proposition}
\begin{proof}
Let $\delta \in \Z_3^\N$. One has
$$\pi^{-1}(\{\delta\})=\{(x_i)_{i\geq 1} \in X : \pi_1(x_1)=\delta \}
=\cup_{x \in \pi_1^{-1}(\{\delta\})} X(x), $$
where $X(x)=\{ (x_i)_{i \geq 1} \in X : x_1=x \}$.
By Lemma \ref{le:preimagen} the set
$X(x)$ is at most countable; $\pi_1^{-1}(\{\delta\})$ is finite by Lemma
\ref{factorpi}.
\end{proof}

Now we construct a countable scrambled set.
Let $n \geq 1$ and $m \geq n$. Remark that
the words $m(n-1)$ and $nn$ can be obtained
as sub-words of $\tau^m_m(0(n-1))$ and $\tau^n_m(0n)$ respectively,
then they are sub-words of $X_m$. So we can define
$$x(n)=\lim_{k\to \infty} \tau_n^k(n).n0n\tau_n(0n)...\tau_n^i(0n)...
\in X_n$$
and
$$y(m,n)=\lim_{k\to \infty} \tau_m^k(m).(n-1)0n\tau_m(0n)...\tau_m^i(0n)...
\in X_m$$
where the symbol after the dot is the zero coordinate.
Clearly $\rho_n(x(n+1))=\rho_n(y(n+1,n+1))=x(n)$ and $\rho_m(y(m+1,n))=y(m,n)$, so
that
$$S=\cup_{n\geq 2}
\{(x(1),...,x(n-1),y(n,n),y(n+1,n),y(n+2,n),...)\}\subset X.$$
\begin{proposition}
The countable set $S$ is a scrambled set of $(X,T)$.
\end{proposition}
\begin{proof}
It is enough to observe that
$(x(n),y(n,n)) \in {\bf A}(X,T)$,
$(x(m),y(m,n)) \in {\bf LY}(X,T)$, $m>n$,
and $(y(m,n),y(m,k)) \in {\bf LY}(X,T)$, $m\geq n$, $m\geq k$,
$m \not = k$.
\end{proof}

\section{A further example}

\noindent Like the example in [Au], p.26, and examples in Section \ref{subst}
the system we describe now is
semi--distal but not almost distal; unlike these examples, it
is an extension of an irrational rotation and has countably many
Li-Yorke pairs.

Let us denote by $\Si$ the circle and consider its representation as the
quotient group $\R/\Z$. Given a real number $t$ we denote by 
$\{t\}$ the fractional part of $t$ (that is, $\{t\}$ is 
the unique representative in $[0,1)$ of the
mod $\Z$ congruence class of $t$). Fix an irrational number 
$\alpha \in [0,1)$  
and consider $R_\alpha: \Si \to \Si$ to be the  
rotation by $\alpha$ on the circle, 
$R_\alpha(t) = \{t + \alpha \}$. 
The system $(\Si,R_\alpha)$ is minimal and isometric, and so is distal.

Let $n_0 = 1$. Choose an increasing sequence of positive integers, 
$\{ n_i : i \geq 1 \}$,
such that the fractional parts
$\{ \{n_j\alpha\} : j\geq 1 \}$
form a decreasing sequence in the unit interval which converges
to $\alpha$.
Observe that $\{n_j\alpha\}$ is the fractional part
of $R_\alpha^{n_j}(0)= R_\alpha^{n_j-1}(\alpha)$.
Now choose  $\beta^1$ so that
$\{n_1\alpha\} <  \beta^1 < 1$ and  $\beta^1$ is not a
rational linear combination of $1$ and $\alpha$.
In particular,  $\beta^1 + \Z$ is not in the orbit of
$\alpha+\Z$ in $\Si$. 
Inductively for $j > 1$,
choose  $\beta^j$ rationally independent of $1$, $\alpha$
and $\beta^k$ for $k < j$ such that
$\{n_j\alpha\}   <  \beta^j < \{n_{j-1}\alpha\}$.
Let $I_j = [\{n_j\alpha\},  \beta^j]$ for $j \geq 1$ and
let $K = \{\alpha \}\cup \bigcup\limits_{j\geq 1} I_j$. Thus, $K$ is a
closed set in the unit interval,
and by projecting we can regard it as a closed set in $\Si$.
The boundary of $K$ consists of 
$\partial_\alpha= \{\{n_j\alpha\} : j \geq 0 \}$
and $\partial_\beta=\{\beta^j : j \geq 1 \}$.
The points of $\partial_\alpha$ lie on the orbit of $\alpha$
while the points in $\partial_\beta$ are all on distinct
orbits different from the $\alpha$ orbit as well.

Put $A = \{ 0, 1 \}$ and let $c_0 : \Si \to A$ be the characteristic
function of $K$. That is,  
$c_0(t) = 1$ if $t \in K$ and $c_0(t)= 0$ otherwise. 
Notice that
$c_0$ is continuous at every point of
$\Si$ except for those on the boundary of $K$.
Let $c_i = c_0 \circ R_\alpha^i$ for $i\in  \Z$. Finally, define
$c:\Si \to A^\Z$ by $c(t)=(c_i(t): i\in \Z)$.
Notice that $c$ defines a commuting map between $(\Si,R_\alpha)$
and $(A^\Z,T)$ where $T$ is the shift map.
Denote by $cont(c)$ the continuity set of $c$. It is given by
$$cont(c)^c=\{\{\alpha + n\alpha\}: n \in \Z \} \cup \bigcup_{j \geq 1} 
\{\{\beta^j+ n \alpha\}: n \in \Z \}.$$
The graph of $c$, $G(c)=\{(t,c(t)): t \in \Si \}$,
is an $R_\alpha \times T$ invariant subset of $\Si \times A^\Z$.
Since $c$ is not continuous, it is not a closed subset. However,
$Gcont(c) = \{ (t, c(t)) : t \in cont(c) \}$
form a dense subset of $G(c)$.

Let $Z$ denote the
closure of $G(c)$ in $\Si \times A^\Z$. This is a $R_\alpha\times T$ 
invariant subset of $\Si \times A^\Z$. We denote by  
$F : Z \to Z$ the  restriction of $R_\alpha \times T$ to $Z$ 
and by $\rho: Z \to \Si$ the coordinate
projection. Clearly, $\rho$ is a factor map from
$(Z, F)$ to $(\Si, R_\alpha)$. For $t \in \Si$ let 
$Z_t=\{ (s,c(s)) \in Z : \rho(s)=t \}$ denote the fiber 
$\rho^{-1}(\{t\})$ in $Z$.

\begin{lemma}\label{ref1}
If $t \in cont(c)$ then $(t, c(t))$ is the unique point of $Z$
such that $\rho(t,c(t))=t$. In particular,
$Z_t$ is a singleton set. More generally, if $(t,x)\in Z$ and $t$
is a continuity point for $c_i$ (that is,  $R_\alpha^i(t)$ is not
in the boundary of $K$) then $x_i = c_i(t)$.
\end{lemma}
\begin{proof}
We leave it as an exercise.
\end{proof}
 
As we observed before the set $Gcont(c)$ is dense in $Z$, then 
$\rho$ is an almost one-to-one map. 
Such maps are always proximal and minimality
is preserved by almost one-to-one lifts. Therefore  $(Z, F)$ is a
minimal system. 
 
For each $t \in cont(c)$ there is a unique $F$ orbit in $Z$
lying over its $R_\alpha$
orbit in $\Si$. The remaining orbits are those of the  $\beta^j$ for $j\geq 1$
and the orbit of $\alpha$.

\begin{lemma}
\label{ref2}
For $j \geq 1$ the fiber $Z_{\beta^j}$ consists of
two points $z^{j+}= (\beta^j, a^j)$ and $z^{j-}= (\beta^j, b^j)$
with $a^j_i=b^j_i=c_i(\beta^j)$ for all $i \geq 1$ and with
$a_0^j=1$ and $b^j_0=0$.
\end{lemma}
\begin{proof}
Let $j \geq 1$.
If a sequence $(t_k)_{k\in \N}$ approaches $\beta^j$ from above then,
eventually, $t_k \notin K$ and so, eventually, $c_0(t_k) = 0$.
For $i \geq 1$ the function $c_i$
is continuous at  $\beta^j$, so the sequence $(c_i(t_k))_{k \in \N}$ 
has limit $c_i(\beta^j)$.
Thus, $((t_k, c(t_k)))_{k \in \N}$ approaches $z^{j-}$. 

If $(t_k)_{k \in \N}$ approaches $\beta^j$
from below then, eventually,
$t_k \in K$, and so $((t_k, c(t_k)))_{k \in \N}$ approaches $z^{j+}$.
Any real sequence admits a monotone subsequence, then these are
the only points in $Z_{\beta^j}$.
\end{proof}

\begin{lemma}\label{ref3}
The fiber $Z_\alpha$ consists of three
points $z^{+1}=(\alpha,a)$, $z^-=(\alpha,b)$ and
$z^{+0}= (\alpha, \bar a)$ such that:
 
$$\text{ for } i\notin \{ n_j-1 : j \geq 0 \},  \ a_i =\bar a_i= b_i =
c_i(\alpha)$$
 
$$ \text{ for } i \in \{ n_j-1 : j > 0 \}, \  a_i =\bar a_i = 1, b_i = 0$$
 
$$ \text{ for } i = 0, \  a_0 = 1, \bar a_0 = b_0 = 0.$$
\end{lemma}
\begin{proof}
If $i \notin\{ n_0 - 1, n_1 -  1, ...\}$ then $c_i$ is continuous at
$\alpha$. If a sequence $(t_k)_{k\in\N}$ approaches $\alpha$ then 
$(c_i(t_k))_{k\in\N}$
approaches $c_i(\alpha)$ for such $i$. 
If $(t_k)_{k\in\N}$ approaches
$\alpha$  from below then it lies in the complement of $K$ and so
$c_0(t_k) = 0$ for all $k\in \N$. For $i = n_j - 1$ with $j > 0$
the fractional part of
$(R_\alpha^{n_j-1}(t_k))_{k\in \NN}$ 
approaches $\{n_j\alpha\}$ from below and so
the rotated sequence is also eventually in the complement of $K$.
Hence $((t_k, c(t_k)))_{k\in\N}$
approaches $z^-$. 
If $(t_k)_{k\in\N}$ approaches $\alpha$ from above then
for $i=n_j-1$ with $j > 0$ the sequence of the fractional parts
of $R_\alpha^{n_j-1}(t_k)$, $k \in \N$,
approaches $\{n_j\alpha\}$
from above and so the rotated sequence eventually enters the interval
$I_j$. Hence, the limit of $(c_i(t_k))_{k\in\N}$ is $1$ for such $i$'s. 
Finally,
for $i = 0$, observe that $\alpha$ is a limit point from above of
points in $K$ and of points in $\Si \setminus K$. If we choose a
sequence $(t_k)_{k\in \N}$ with $t_k \in K$ for
all $k\in \N$ then $c_0(t_k) = 1$ and the 
sequence $((t_k, c(t_k)))_{k\in \N}$ approaches $z^{+1}$.
If instead we choose a sequence $(t_k)_{k\in \N}$ with $t_k \notin K$ 
for all $k\in \N$ then
$c_0(t_k)=0$ and $((t_k, c(t_k)))_{k\in\N}$ approaches $z^{+0}$.
Any sequence approaching
$\alpha$ has a subsequence of one of these three types. Hence,
the fiber consists of the three points specified in the statement.
\end{proof}

Since $(\Si, R_\alpha)$ is a distal system, a proximal pair of distinct points
in $Z$ must lie in the same fiber of $\rho$.
 
It is clear that the pair of orbits which lie over the orbit of
$\beta^j$ is an asymptotic pair for each $j\geq 1$. 
There remain the three
orbits which lie over the orbit of $\alpha$. Obviously, the pair
$(z^{+0}, z^{+1})$ is an asymptotic pair, while the pairs
$(z^{+0}, z^-)$
and $(z^{+1}, z^-)$ are proximal pairs, but not asymptotic.
Thus, the Li-Yorke set is the countable collection consisting
of the orbits of these two pairs.

Now by Akin's remark (see the proof of Corollary \ref{coroequiv})
and the last discussion 
we conclude that $(Z,T \times \sigma)$
contains no strong Li--Yorke pairs.

\medskip
{\bf Acknowledgments.} We thank Pierre Arnoux and
Xiangdong Ye for fruitful discussions and comments.
We also thank all the valuable comments and suggestions 
of the anonymous referees which in 
particular help to give the final form  
to the example in Section 5.
The authors acknowledge financial support
from Nucleus Millennium Information and Randomness P01-005,
FONDECYT 1010447 and ECOS-Conicyt cooperation agreement C99E10.

\end{document}